\documentclass[12pt]{article}

\usepackage{amssymb}
\usepackage{amscd}
\usepackage{amsthm}
\usepackage{amsmath}
\usepackage{latexsym}
\usepackage{graphicx}
\usepackage{delarray}
\usepackage{epic}

\begin{document}


\newtheorem{theorem}{Theorem}[section]
\newtheorem{definition}[theorem]{Definition}
\newtheorem{lemma}[theorem]{Lemma}
\newtheorem{proposition}[theorem]{Proposition}
\newtheorem{corollary}[theorem]{Corollary}
\newtheorem{example}[theorem]{Example}
\newtheorem{remark}[theorem]{Remark}

\newcommand{\G}{\mathcal{G}}
\newcommand{\FT}{\mathcal{F}}
\newcommand{\Hilb}{\mathcal{H}}
\newcommand{\dG}{{\hat{\mathcal{G}}}}
\newcommand{\g}{\gamma}
\newcommand{\dg}{d\gamma}
\newcommand{\Lam}{\Lambda}
\newcommand{\lam}{\lambda}
\newcommand{\dLam}{\Lambda^\perp}
\newcommand{\R}{\mathbb{R}}
\newcommand{\N}{\mathbb{N}}
\newcommand{\Z}{\mathbb{Z}}

\newcommand{\set}[2]{\big\{ \, #1 \, \big| \, #2 \, \big\}}
\newcommand{\inner}[2]{\langle #1,#2\rangle}
\def\supp{\mbox{\rm supp}}


\title{Minimal Norm Interpolation in Harmonic Hilbert Spaces and
Wiener Amalgam Spaces on Locally Compact Abelian Groups}
\author{
 H.~G.~Feichtinger
 \footnote{NuHAG, Department of Mathematics,
 University of Vienna, Nordbergstrasse 15, \newline A-1090 Vienna,
 Austria. \{hans.feichtinger, tobias.werther\}@univie.ac.at \newline
The third author gratefully acknowledges partial support by the
FWF-project 15605 of the Austrian Science Foundation.} \, ,
 S.~S.~Pandey
 \footnote{Department of Mathematics, R.~D.~University,
 Jabalpur, India. \newline sheelpandey@hotmail.com}\, , and
 T.~Werther $^*$
}


\maketitle



\vspace{.2cm}



\begin{abstract}
The family of harmonic Hilbert spaces is a natural enlargement of
those classical $L^{2}$-Sobolev space on $\R^d$ which consist of
continuous functions. In the present paper we demonstrate that the
use of basic results from the theory of Wiener amalgam spaces
allows to establish fundamental properties of harmonic Hilbert
spaces even if they are defined over an arbitrary locally compact
abelian group $\G$. Even for $\G = \R^{d}$ this new approach
improves previously known results. In this paper we present
results on minimal norm interpolators over lattices and show that
the infinite minimal norm interpolations are the limits of finite
minimal norm interpolations. In addition, the new approach paves
the way for the study of stability problems and error analysis for
norm interpolations in harmonic Hilbert and Banach spaces on
locally compact abelian groups.
\end{abstract}


\section{Introduction}

Babuska \cite{Bab68} introduced the concept of periodic Hilbert
spaces in order to study universally optimal quadrature formulas.
Subsequently Prager \cite{Pr79} has studied in detail the
relationship between optimal approximation of linear functionals
on periodic Hilbert spaces and minimal norm interpolation. In a
more recent paper Delvos \cite{Del97} has introduced the notion of
harmonic Hilbert spaces for the real line $\R$ and discussed the
interpolation problems over an infinite uniform lattice for
$\ell^2$-data in $\R$.

In the present paper we introduce the concept of harmonic Hilbert
spaces over a locally compact abelian (LCA) group $\G$, establish
the properties of these spaces, and provide a characterization of
minimal norm interpolators. The use of Wiener Amalgam spaces
provides a clear understanding of the subject. This paper opens
new viewpoints for the study of optimal norm interpolation in
various types of harmonic Banach spaces over locally compact
abelian groups.

The article is organized as follows. In Section \ref{sec:pre} we
introduce the notation and state basic facts of harmonic analysis
on compact abelian groups that we essentially need for proving the
main statements. The next section gives a short overview on Wiener
Amalgam spaces and contains a summary of some of their important
properties. In Section \ref{sec:rkhs} we introduce weight
functions and derive the reproducing kernel of harmonic Hilbert
spaces. Then we show in Section \ref{sec:riesz} how the concept of
Riesz-basis leads to certain principal shift invariant subspaces.
In particular, it follows that the bi-orthogonal system of the
Riesz sequence of uniform translates of the kernel is generated by
the same family of translates of the so-called Lagrange
interpolator. Section \ref{sec:int} contains the core results of
the paper, namely the explicit derivation of minimal norm
interpolations. In the final section we discuss the approximation
of the minimal norm interpolation by optimal interpolations of
finite samples.


\section{Preliminaries and Notation}\label{sec:pre}

Let $\G$ be a locally compact abelian group and $\dG$ its dual
group with the normalized Haar measures $dx$ and $d\g$,
respectively. The dual group $\dG$ is defined as the set of all
characters $\g:\G\to \mathrm{T}$, i.e., continuous homomorphisms
from $\G$ into the torus $\mathrm{T}$ with pointwise
multiplication as group operation, and the standard (compact-open)
topology.

For $1\leq p\leq \infty$, $L^p$-spaces are denoted by their usual
symbols. For $f\in L^1(\G)$, the Fourier transform $\FT$ is
defined as
 \[
 \FT f(\g) = \hat{f}(\g) = \int_\G f(x)\overline{\g(x)} dx\,.
 \]
There are many excellent textbooks on abstract harmonic analysis
such as \cite{RS00} and \cite{Fol95} that comprise a rigorous
treatise of Fourier transform on LCA groups. Here, we only state
those facts that are most useful for the development of the
presented results.

Using the Pontryagin duality theorem which allows to write
$\inner{\g}{x}$ or $\inner{x}{\g}$ unambiguously for $\g(x)$, the
Fourier inversion formula for $h=\hat{f}\in L^1(\dG)$ writes as
\[ f(x) = \FT^{-1}h(x) = \int_\dG \hat{f}(\g)
\inner{\g}{x}d\g\,.\]
 The range of the Fourier transform, which is an injective and
bounded linear mapping from $L^1$ to $C^0$, the space of
continuous functions that vanish at infinity, is denoted by $\FT
L^1$. The convolution theorem states that \[ \FT(f\ast g) =
\FT(f)\FT(g),\qquad f,g\in L^1\,, \] where $"\ast"$ denotes the
usual convolution. By transport of the norm $\|\FT f\|_{\FT L^1} =
\|f\|_{L^1}$, the normed space $(\FT L^1,\|\cdot\|_{\FT L^1})$
turns into a Banach algebra with respect to pointwise
multiplication. Due to Plancherel's theorem one has $\FT L^1 = L^2
\ast L^2$, which ensures that $\FT L^1$ is a dense subspace of
$C^0$ endowed with the sup-norm $\|\cdot\|_\infty$.

Whenever a Banach space $(X,\|\cdot\|_X)$ is continuously embedded
into another Banach space $(Y,\|\cdot\|_Y)$, i.e., there exists a
constant $C>0$ such that $\|x\|_Y\leq C\|x\|_X$ for all $x\in X$,
we write $X\hookrightarrow Y$.

Throughout this paper we use the symbol $T_{x}$ for the
translation operator \[ T_{x} f(y) = f(y-x), \quad x,y \in\G. \]
For convenience, all index sets, sums, and lattices on $\G$ and
$\dG$ used in the sequel are assumed to be countable.

Let $\Lam$ be a subgroup of $\G$. According to \cite{Gro98}, we
call $\Lam$ a lattice if the quotient group $\G/\Lam$ is compact.
The lattice size $s(\Lam)$ is defined as the measure of a
fundamental domain of $\Lam$ in $\G$, i.e., we choose a measurable
set $U\subset \G$ such that every $x\in\G$ can be uniquely written
as $x=\lam+u$ for some $\lam\in\Lam$ and $u\in U$,
cf.~\cite{Gro98}. Then the lattice size $s(\Lam)$ is the measure
of $U$. For an equivalence class in $\G/\Lam$ we write $\dot{x} =
\{x+\lam\}$. Let $d\dot{x}$ denote the normalized Haar measure on
$\G/\Lam$. For $f\in L^1(\G)$, we then have Weil's formula
\cite{RS00}
\begin{equation}\label{eq:weil1}
\int_\G f(x)dx = s(\Lam)\int_{\G/\Lam}
\Big(\sum_{\lam\in\Lam}f(x+\lam) \Big)d\dot{x}\,.
\end{equation}

The annihilator of $\Lam$ is the subgroup \[ \Lam^\perp =
\set{\chi\in\dG}{\chi(\lam)=1,\; \lam\in\Lam} \] and Weil's
formula becomes in this context the identity
\begin{equation}\label{eq:weil2}
\int_\dG \hat{f}(\g)d\g = s(\Lam^\perp)\int_{\dG/\Lam^\perp}
\Big(\sum_{\chi\in\Lam^\perp}\hat{f}(\g+\chi) \Big)d\dot{\g}\,.
\end{equation}

Note that the dual group $\hat{\Lam}$ of $\Lam$ is naturally
isomorphic to $\dG/\Lam^\perp$, \cite{RS00}. For
$(c_\lam)\in\ell^2(\Lam)$, the sum $\sum_{\lam\in\Lam}c_\lam
\inner{\g}{\lam}$ is a Fourier series on $\dG/\Lam^\perp$ and
Plancherel's theorem \cite{RS00} yields
\begin{equation}\label{eq:plancherel}
\int_{\dG/\Lam^\perp}\Big|\sum_{\lam\in\Lam}c_\lam
\inner{\g}{\lam}\Big|^2 d\dot{\g} = \sum_{\lam\in\Lam}|c_\lam|^2
\,.
\end{equation}
For $\G=\R$ and $\Lam=\alpha\Z$, $\alpha > 0$, this is just the
classical Fourier series expansion of periodic functions. In this
case we simply have $\Lam^\perp = \alpha^{-1}\Z$.

We make use of these fundamental results in the proof of Theorem
\ref{th:riesz}.


\section{Wiener Amalgam Spaces}

Wiener Amalgam spaces on LCA groups and their properties under
Fourier transform are studied in a series of papers starting with
\cite{Fei80}. We refer to \cite{ Fei90, Hei03} for a compact
survey on their properties. Here, we only provide a short overview
over the facts needed in the present context.

First we introduce the continuous version of Wiener Amalgam
spaces. Let $\psi$ be a non-zero, non-negative function on $\G$
with compact support, i.e., $\supp(\psi)\subset \Omega$ where
$\Omega \subseteq G$ is compact. Given any Banach space $(B, \|
\cdot \|_B)$ of functions on $\G$ on which $\FT L^1$ acts by
pointwise multiplication, i.e., $\FT L^1 \cdot B \subseteq B$, the
Wiener Amalgam space $W(B,L^q)$ is defined by
\begin{equation}\label{eq:wam1b}
 W(B,L^q) = \set{f\in B_{loc}}
 { \|f\|_{W(B,L^q)} := \Big ( \int_\G  \| f
 \cdot T_x \psi  \|_B^{q} dx \Big )^{1/q}  <\infty}
\end{equation}
for $1\leq q < \infty$ (with obvious modifications for
$q=\infty$). The expression $f\in B_{loc}$ means that the function
$f$ can at least locally be measured by the $B$-norm hence the
integrand in the definition of $\|\cdot\|_{W(B,L^q)}$ is a
well-defined non-negative function. As stated in \cite{Fei90},
these spaces are Banach spaces and do not depend on $\psi$ (up to
equivalence of norms).

An equivalent, but discrete, definition of Wiener Amalgam spaces
uses a so called $\FT L^1$-bounded, uniform partition of unity
(for short BUPU), that is a sequence of non-negative functions
$\{\psi_i\}$ corresponding to a sequence $\{g_i\}$ in $\G$ such
that
\begin{enumerate}
 \item $\supp(\psi_i) \subset g_i + \Omega$,
 \item $\sup_i\sharp\set{j}{(g_i+\Omega)\cap (g_j+\Omega) }<\infty$,
 \item $\sum_i \psi_i(g) \equiv 1$,
 \item $\psi_i$ are bounded in $\FT L^1$, i.e.,
       $\|\psi_i\|_{\FT L^1}\leq C <\infty$ for all $i$.
\end{enumerate}
With the help of such a BUPU we define the discrete Wiener Amalgam
space
\begin{equation}\label{eq:wam2}
W(B,\ell^q) = \set{f \in B_{loc} }{\|f\|_{W(B,\ell^q)} := \Big(
\sum_i \|f \psi_i\|_B^q \Big)^{1/q} < \infty }
\end{equation}
for $1\leq q < \infty$. Note that again, $W(B,\ell^q)$ is a Banach
space independent from the partition of unity (up to equivalence
of norms). In all cases of this paper we will use a lattice
$\Lambda$ in $\G$ instead of an arbitrary sequence $\{g_i\}$.

In contrast to the $L^p$-spaces, where, for instance, on $\R^d$,
inclusion fails for different $p$, Wiener Amalgam spaces enjoy the
so-called {\it coordinate-wise inclusion}, i.e., if
$B_{1,loc}\hookrightarrow B_{2,loc}$ and $q_1\leq q_2$, then
$W(B_1,\ell^{q_1})\hookrightarrow W(B_2,\ell^{q_2})$.

Wiener Amalgam spaces can also be defined with additional weight
functions that are described in the following section,
cf.~\cite{Fei90}.


\section{Reproducing Kernel Hilbert Spaces}\label{sec:rkhs}

We first state sufficient conditions on weight functions that are
used to define harmonic Hilbert spaces. It is beyond the scope of
this work to elaborate on these conditions.

A strictly positive and continuous function $w$ is called a {\it
submultiplicative}  (or {\it Beurling}) weight function on $\dG$
if
\begin{equation}\label{eq:weight}
 w(\g_1+\g_2) \; \leq  \; w(\g_1)\, w(\g_2) \quad  \mbox
{for all} \quad \g_1,\g_2 \in \dG,
\end{equation}
It satisfies the {\it Beurling-Domar non-quasianalyticity condition}
\cite[VI, sect.~3]{RS00} if
\begin{equation}\label{eq:beurldomar}
\begin{array}{ll}
  \mbox{(BD)}& \qquad \overset{\infty}{\underset{n=1}{\sum}}
  n^{-2} \log w(n\g) < \infty
\quad \mbox {for all} \quad \g\in\dG \,. \\
\end{array}
\end{equation}
The standard example of such weight functions are weights of
polynomial type on $\R^d$ such as
\[ w_s(\g) = (1+|\g|)^s \simeq (1+|\g|^2)^{s/2}\,, \qquad
\g\in\R^d, \] for $s\geq 0$, or subexponential weights such as
\[ w(\g) = e^{\alpha |\g|^\delta}(1+|\g|)^s
\,, \qquad \g\in\R^d, \] for $\alpha>0,\; 0<\delta<1$, and $s\in\R$. For
detailed studies on such weight functions on the Euclidian space $\R^d$ we
refer to \cite{Fei79,Fei83,FeiMod03}.

Throughout this paper, we assume that $w^{-1} \in L^2(\dG)$ which implies
$w^{-2}\in W(C^0,\ell^1)$. Obviously, this is satisfied for the above example
$w_s$ whenever $s>d/2$. In contrast to \cite{Del97} where the author included
the box function as a possible weight leading to band-limited functions that
we discuss separatly in Example 2, we restrict our discussion to
submultiplicative weight functions satisfying (BD) in order to obtain
stronger results and general statements valid for arbitrary lattices. From
our point of view the band-limited case should be seen as a limiting case,
requiring sometimes separate arguments.

Among others, the submultiplicativity of $w$ in conjunction with the
assumption  $w^{-2} \in L^1(\dG)$, allows to apply the convolution theorem
for Wiener Amalgam spaces \cite{Fei80,Hei03}, in order to derive the
following crucial property: for any lattice $\Lam$ there exist positive
constants $a$ and $b$ such that
\begin{equation}\label{eq:periodize}
a\leq \sum_{\chi\in\Lam^\perp}w^{-2}(\g+\chi)\leq b\,, \qquad
\g\in\dG\,.
\end{equation}
The upper bound follows from the fact that the Haar measure for
$\Lambda^\perp$, i.e.\ $ \mu = \sum_{\chi\in\Lam^\perp}
\delta_\chi$ belongs to $W(M,\ell^\infty)$, where $M$ denotes the
space of bounded measures, while on the other hand the
submultiplicativity implies that $w^{-2} \in W(C^0,\ell^1)$. Since
the $\Lambda^{\perp}$ periodization of $w^{-2}$ equals $\mu \ast
w^{-2} \subseteq W(M,\ell^\infty) \ast W(C^0,\ell^1) \subseteq
W(C,\ell^\infty) = C^b(G)$ the upper bound is valid, see
\cite{Fei90}. Since $w^{-2}$ is strictly positive the periodicity
implies that there is a strictly positive lower bound.

We define the {\it harmonic Hilbert} space corresponding to the
weight $w$ as
\begin{equation}\label{eq:hh}
\Hilb_w(\G) =  \FT^{-1}L^2_w(\dG)
\end{equation}
with inner product
 \[
 \inner{f}{g}_w = \int_\dG \hat{f}(\g)
 \overline{\hat{g}(\g)}w^2(\g)d\g\,.
 \]
\begin{lemma}
$\Hilb_w(\G)$ is a Hilbert space.
\end{lemma}
\begin{proof}
Left to the reader.
\end{proof}
In the case of the weight functions $w_s$, $\Hilb_{w_s}$ is known
as the Sobolev space of fractional order $s$ \cite{Ada75}.

\noindent{\bf Remark.} It has been shown in \cite{Fei79} that if
$w^{-1}\in L^2$ is {\it subadditive}, i.e.,
\[ w(\g_1+\g_2) \leq C(w(\g_1)+w(\g_2))\,, \qquad
\g_1,\g_2\in\dG\,,\] then $L^2_w$ is a Banach convolution algebra.
As a consequence, $\Hilb_w$ turns into a Banach algebra with
respect to pointwise multiplication. This is the case for $w_s$,
$s>d/2$, and subexponential weights.

We now state a result which is fundamental for sampling functions
in $\Hilb_w$ along some lattice $\Lam$.
\begin{theorem}\label{th:embedding}
If $w^{-1} \in L^2(\dG)$, then
 \[\Hilb_w(\G)  = W(\FT^{-1}L^2_w,\ell^2)
 \hookrightarrow W(C^0, \ell^2)(\G) \,. \]
\end{theorem}
\begin{proof}
By virtue of Corollary 7 of \cite{Fei90}, we have $\FT^{-1}L^2_w =
W(\FT^{-1}L^2_w,\ell^2)$. We emphasize that this result is a consequence of
the submultiplicativity and (BD). Now, since $w^{-1}\in L^2$, an easy
application of the Cauchy-Schwartz inequality implies  $L^2_w \hookrightarrow
L^1$ which in turn results in the following inclusions
 \[
\FT^{-1}L^2_w \hookrightarrow \FT^{-1}L^1 \hookrightarrow C^0,
 \]
(by the Riemann-Lebesgue Lemma). Hence, by the coordinate-wise inclusion
properties of Wiener Amalgam spaces we obtain
\[\Hilb_w = W(\Hilb_w,\ell^2) \hookrightarrow
W(C^0, \ell^2)\,.\]
\end{proof}
As a consequence, $\Hilb_w(\G)$ is a {\it reproducing kernel
Hilbert space} (RKHS). This is a Hilbert space consisting of
continuous functions in which the point evaluation functionals are
continuous. Therefore, by virtue of the Riesz Representation
theorem, for each $x\in\G$ there exists a unique function, say
$k_x$, in the RKHS, such that the point evaluation at $x$ of any
function $f$ in the RKHS can be performed by means of the inner
product with $k_x$. In the case of $\Hilb_w$, that is
 \[
f(x)= \inner{f}{k_x}_w\,, \qquad f\in\Hilb_w(\G)\,, x\in\G\,.
 \]
The kernel $k(x,y) = k_x(y)$  defines a continuous function on $\G \times\G$
containing all the informations about the scalar product. A first detailed
survey on RKHS goes back to \cite{Aro50}.

In the present situation the kernel of $\Hilb_w$ turns out to consist of
translations of a single function.

\begin{proposition}\label{cor:rkh}
$\Hilb_w(\G)$ is a reproducing kernel Hilbert space with kernel
$k(x,y) = \phi(x-y)$, where \[ \phi = \mathcal{F}^{-1} w^{-2}\,.
\] That is, $f(x) = \inner{f}{T_x\phi}_w$ for all
$f\in\Hilb_w$ and $x\in\G\,. $
\end{proposition}
\begin{proof}
Define $\phi$ by $\hat{\phi}=w^{-2}$. Since $w^{-2}\in L^1$, \,
$\phi\in\Hilb_w$\,~. Next we compute
\begin{eqnarray}
\inner{f}{T_x\phi}_w & = & \int_\dG \hat{f}(\g)
\overline{\widehat{T_x\phi}(\g)}w^2(\g)d\g \nonumber \\
&=& \int_\dG \hat{f}(\g)
\overline{\overline{\inner{x}{\g}}w^{-2}(\g)}w^2(\g)d\g \label{eq:eval}\\
&=& \int_\dG \hat{f}(\g)\inner{x}{\g}d\g  \;= \;f(x)\,, \nonumber
\end{eqnarray}
the last step following from the Fourier inversion theorem (note that $\hat f
\in L^1$). The uniqueness of the kernel completes our proof of the fact $k_x
= T_x \phi$.
\end{proof}

\noindent{\bf Remark.} The fact that $\Hilb_w\hookrightarrow C^0$ can also be
derived immediately  from (\ref{eq:eval}). The use of Wiener Amalgam spaces,
however, reveals the important property that the sequence of samples of any
function in $\Hilb_w$ with respect to any lattice $\Lam$ is square summable.
This property has led to the definition of so-called {\it $\ell^2$- puzzles}
in \cite{Tch84}.

\begin{corollary}\label{cor:map}
For any lattice $\Lam$ the mapping
\begin{equation}\label{eq:mapQ}
Q:f\mapsto (f(\lam))_{\lam\in\Lam}
\end{equation}
is bounded from $\Hilb_w(\G)$ to $\ell^2(\Lam)$.
\end{corollary}
\begin{proof}
Choose a fixed lattice $\Lam'$ of $\G$ with fundamental domain
$\Omega'$. Then the family $\{T_{\lam'}
\chi_\Omega'\}_{\lam'\in\Lam'}$ where $\chi$ denotes the
characteristic function, forms a BUPU. We now consider the norm of
$W(C,\ell^2)$ with respect to this BUPU. For any $\lam'\in\Lam'$,
there is at most a finite number of points of $\Lam$ in
$\lam'+\Omega'$, say $n_{\lam'}$. These numbers are uniformly
bounded. Therefore, $\|(f(\lambda))\|_{\ell^2} \leq
\sup(n_{\lam'})\|f\|_{W(C,\ell^2)}$. Hence, the result follows
from Theorem \ref{th:embedding}.
\end{proof}

We will later see that the mapping $Q$ is surjective.


\section{Riesz Basis}\label{sec:riesz}

Riesz bases are a well-established concept in Hilbert space theory
\cite{Chr03}.
\begin{definition}\label{def:riesz}
A family of vectors $\{h_n \}$ in a Hilbert space $\Hilb$ is
called a Riesz sequence if there exist bounds $0<a\leq b<\infty$
such that
\begin{equation}\label{eq:Riesz-bounds}
a\|c\|^2_{\ell^2}\leq\big\|\sum_n c_n h_n  \big\|^2\leq
b\|c\|^2_{\ell^2}
\end{equation}
for all sequences $c=(c_n)\in\ell^2$.
\end{definition}
Riesz sequences generalize the concept of orthogonal sequences as
one can see from the following properties. Let us call $V$ to be
the closed linear span of $\{h_n\}$. For every Riesz sequence
$\{h_n \}$ there exists a unique {\it dual} or {\it bi-orthogonal
sequence} $\{\tilde{h}_n \}$ in $V$ such that
\[ \inner{h_n}{\tilde{h}_m} = \delta_{nm}\,, \qquad n,m\in\N\,, \]
and the orthogonal projection $P_V$ from $\Hilb$ onto $V$ is given by
\begin{equation}\label{eq:riesz-sum}
P_V h = \sum_n\inner{h}{\tilde{h}_n}h_n =
\sum_n\inner{h}{h_n}\tilde{h}_n \,,\quad\mbox{for all}\quad
h\in\Hilb\,,
\end{equation}
cf.~\cite{You01}. A Riesz sequence $\{h_n \}$ obviously
constitutes a (Riesz) basis for $V$.

The following statement is a standard result in Fourier analysis,
cf.~\cite{RS95}. For the sake of completeness we include the
prove.
\begin{theorem}\label{th:riesz}
For any lattice $\Lam$ of $\G$, the sequence $\{T_\lam
\phi\}_{\lam\in\Lam}$ is a Riesz sequence in $L^2(\G)$ if and only
if there exist positive constants $a,b$ such that
\begin{equation}\label{eq:period1}
a\leq \sum_{\chi\in\Lam^\perp}|\hat{\phi}(\g+\chi)|^2\leq
b\qquad\mbox{a.~e.}
\end{equation}
\end{theorem}
\begin{proof}
We compute
\begin{eqnarray*}
&&\big\| \sum_{\lam\in\Lam} c_\lam T_\lam\phi \big\|^2_2
\stackrel{\rm Plancherel}{=} \big\| \sum_{\lam\in\Lam} c_\lam
\g_\lam \hat{\phi} \big\|^2_2 \\ &=& \int_\dG
\big|\hat{\phi}(\g)\big|^2\big|\sum_{\lam\in\Lam} c_\lam
\inner{\lam}{\g} \big|^2 d\g \\ &\stackrel{\rm
(\ref{eq:weil2})}{=}&
s(\Lam^\perp)\int_{\dG/\Lam^\perp}\sum_{\chi\in\Lam^\perp}\Big[
\big|\hat{\phi}(\g+\chi)\big|^2\big|\sum_{\lam\in\Lam} c_\lam
\inner{\lam}{\g+\chi} \big|^2\Big]d\dot{\g} \\ &=&
s(\Lam^\perp)\int_{\dG/\Lam^\perp}\sum_{\chi\in\Lam^\perp}\Big[
\big|\hat{\phi}(\g+\chi)\big|^2\big|\sum_{\lam\in\Lam} c_\lam
\inner{\lam}{\g}\inner{\chi}{\lam} \big|^2\Big]d\dot{\g} \\ &=&
s(\Lam^\perp)\int_{\dG/\Lam^\perp}\sum_{\chi\in\Lam^\perp}\Big[
\big|\hat{\phi}(\g+\chi)\big|^2\big|\sum_{\lam\in\Lam} c_\lam
\inner{\lam}{\g}\big|^2\Big]d\dot{\g} \\ &=&
s(\Lam^\perp)\int_{\dG/\Lam^\perp}\sum_{\chi\in\Lam^\perp}
\big|\hat{\phi}(\g+\chi)\big|^2\big|\sum_{\lam\in\Lam} c_\lam
\inner{\lam}{\g}\big|^2d\dot{\g}\,.
\end{eqnarray*}
Since this holds for all $(c)_\lam\in\ell^2(\Lam)$, the statement
follows from (\ref{eq:plancherel}).
\end{proof}
A similar result holds for $\Hilb_w(\G)$ instead of $L^2(\G)$.
\begin{theorem}\label{th:riesz2}
For any lattice $\Lam$ of $\G$, the sequence $\{T_\lam
\phi\}_{\lam\in\Lam}$ is a Riesz sequence in $\Hilb_w(\G)$ if and
only if there exist positive constants $a,b$ such that
\begin{equation}\label{eq:period2}
a\leq \sum_{\chi\in\Lam^\perp}\hat{\phi}(\g+\chi)\leq
b\qquad\mbox{a.~e.}
\end{equation}
\end{theorem}
\begin{proof} Analogue to the proof of Theorem \ref{th:riesz}.
\end{proof}

Recalling Condition (\ref{eq:periodize}) on the weight, we finally
obtain that $\{T_\lam \phi\}_{\lam\in\Lam}$ forms a Riesz basis
for its closed  linear span, the so-called {\it spline-type space}
 \[
V_\Lam(\phi) \; = \; \mbox{span}\set{T_\lam \phi}{\lam\in\Lam} \;
= \; \set{ f = \sum_{ \lam\in\Lam}c_\lam T_\lam\phi}
{(c_\lam)\in\ell^2(\Lam)}
 \]
which is a closed subspace of $\Hilb_w(\G)$.

We now look for the dual basis of $\{T_\lam \phi\}_{\lam\in\Lam}$
in $\Hilb_w$. Let us define
\begin{equation}\label{eq:dualatom}
\hat{\psi}(\g) =
\frac{\hat{\phi}(\g)}{\sum_{\chi\in\Lam^\perp}\hat
{\phi}(\g+\chi)}\,.
\end{equation}
Note that the denominator is a $\Lam^\perp$-periodic square
integrable function. Since any such function has a unique Fourier
series presentation with $\ell^2$-coefficients, we can easily see
that $\psi$ belongs to $V_\Lam(\phi)$. Similar to the proof of
Theorem \ref{th:riesz} we compute
\begin{equation}\label{eq:lagrprop}
\inner{T_\lam \phi}{T_{\lam'}\psi}_w =
\delta_{\lam,\lam'}\,,\qquad\lam,\lam'\in\Lam\,.
\end{equation}
Hence, $\{T_\lam \psi\}_{\lam\in\Lam}$ is the dual basis of
$\{T_\lam \phi\}_{\lam\in\Lam}$ in $V_\Lam(\phi)$. It follows from
(\ref{eq:riesz-sum}) that every function in $V_\Lam(\phi)$ can be
written as
\begin{equation}\label{eq:expansion}
f = \sum_{\lam\in\Lam}\inner{f}{T_\lam\phi}\, T_\lam \psi =
\sum_{\lam\in\Lam}f(\lam)\, T_\lam \psi\,.
\end{equation}
As a consequence, every function in $V_\Lam(\phi)$ is completely
determined by its samples on $\Lam$.

Because of the reproducing kernel property of $\phi$,
(\ref{eq:lagrprop}) is equivalent to
\[
 \psi(\lam-\lam') = \delta_{\lam,\lam'}\,,\qquad \lam,\lam'\in\Lam\,.
\]
Therefore, $\psi$ is also called the {\it Lagrange interpolator}
for $\Lam$.

It is important to note that the bi-orthogonal system of $\{T_\lam
\phi\}_{\lam\in\Lam}$ is again generated by the same translates of
a single function, namely the Lagrange interpolator.

\noindent{\bf Remark.} Following Theorem \ref{th:riesz},
$V_\Lam(\phi)$ is also a closed subspace of $L^2$. Similar to
above, the dual basis of $\{T_\lam \psi\}_{\lam\in\Lam}$ in $L^2$
is generated by the dual atom $\psi_2$ given by
 \[
\hat{\psi_2}(\g) =
\frac{\hat{\phi}(\g)}{\sum_{\chi\in\Lam^\perp}|\hat
{\phi}(\g+\chi)|^2}\,.
 \]
\noindent{\bf Example 1.} In \cite{FW02}, the authors give a
detailed study of interpolation and stable reconstruction of
functions in the Sobolev space $\Hilb_{w_s}(\R^d)$ from samples
taken over a lattice of varying lattice size.

\noindent{\bf Example 2.} A well-known example is the space of
band-limited functions \[ B = \set{f\in L^2(\R)}
{\mbox{supp}(\hat{f}) \subset [-1/2,1/2]}\] endowed with the
$L^2$-inner product. It might be seen as a Harmonic Hilbert space
for the weight function
 \[
w(x) = \left\{
\begin{array}{cc}
  1 & x\in [-1/2,1/2]\,,  \vspace{.2cm} \\
  \infty & x \not\in [-1/2,1/2] \\
\end{array}
\right.
 \]
although the weight function does not satisfy the conditions
(\ref{eq:weight}) and (\ref{eq:beurldomar}). Nevertheless, results
similar to those above can be obtained. For instance, $B$ is a
RKHS with the so-called sinc-kernel
 \[
\mbox{sinc}(x) = \frac{\sin \pi x}{\pi x}
 \]
which is the inverse Fourier transform of the box function
 \[
\xi(x) = \left\{
\begin{array}{cc}
  1 & x\in [-1/2,1/2]\,,  \vspace{.2cm} \\
  0 & x \not\in [-1/2,1/2]\,. \\
\end{array}
\right.
 \]
Note that $\xi = w^{-2}$. The set of integer shifts of the
sinc-kernel constitutes an orthonormal sequence. It is even an
orthonormal basis of $B$. Consider the lattice $\Lam=\alpha\Z$
with $\alpha<1$. The corresponding annihilator is $\Lam^\perp =
\alpha^{-1}\Z$, and it can easily be seen that \[\sum_{k\in\Z}
w^{-2}(x-\alpha^{-1}k)\] is not bounded away from zero. In
particular $\{T_{\alpha k}\mbox{sinc}\}_{k\in\Z}$ is an
overcomplete basis system, a so-called {\it frame} \cite{Chr03}
and not a Riesz basis, in accordance to Theorem \ref{th:riesz}.
The fact that the periodization that is not bounded away from zero
almost everywhere, leads to a frame, is a general property and
characterization of such frames as shown in \cite{BL98}.


\section{Minimal Norm Interpolation}\label{sec:int}

The problem of minimal norm interpolation in a harmonic Hilbert
space on the real line has been discussed by Delvos, \cite{Del97}.
In the present section we study the corresponding problem in a
more general setting by integrating results stated above.

\begin{theorem}\label{th:minp}
Let $\Lam$ be a lattice in $\G$ and $(c_\lam)\in\ell^2(\Lam)$. The
interpolation problem $f(\lam) = c_\lam$ for all $\lam\in\Lam$ has
a unique minimal norm solution in $\Hilb_w(\G)$. It coincides with
the interpolating element in $V_\Lam(\phi)$ with $\hat{\phi} =
w^{-2}$.
\end{theorem}
\begin{proof}
Since Theorem \ref{th:embedding} shows that $\Hilb_w(\G)
\hookrightarrow W(C, \ell^2)(\G)$, the mapping
\[f\in\Hilb_w(\G) \to (f(\lam))\in\ell^2(\Lam)\] is well-defined
and  bounded, cf.~Corollary \ref{cor:map}. Due to Theorem
\ref{th:riesz2}, $\{T_\lam\phi\}$ is a Riesz basis for
$V_\Lam(\phi)$ which, by duality, consists of all functions of the
form \[ f = \sum_{\lam\in\Lam}c_\lam T_\lam \psi\] for some
$(c_\lam)\in\ell^2$ and the Lagrange interpolator $\psi$ which is
the dual Riesz atom for $\{T_\lam\phi\}$. Due to the Lagrange
property of $\psi$ it follows that
\[ f(\lam) = c_\lam\,,\qquad \lam\in\Lam\,, \] for any such $f\in
V_\Lam(\phi)$. Hence the mapping \[f\in V_\Lam(\phi) \to
(f(\lam))\in\ell^2(\Lam)\] is bijective. In particular, there
exists a unique element $f_c\in V_\Lam(\phi)$ with
$f(\lam)=c_\lam$. Since $V_\Lam(\phi)$ is closed, we can split
$\Hilb_w(\G)$ into the direct sum
\begin{equation}\label{eq:direct_sum}
\Hilb_w(\G) = V_\Lam(\phi)\oplus V^\perp_\Lam(\phi)\,,
\end{equation}
where $V^\perp_\Lam(\phi) = \set{f\in\Hilb_w(\G)}{f(\lam)=0,\;
\lam\in\Lam}\,.$ Assume that some $g\in\Hilb_w(\G)$ interpolates
$c$ on $\Lam$. Then $g-f_c\in V^\perp_\Lam(\phi)$ and we obtain
\[\|g\|^2_w = \|(g-f_c) + f_c\|^2_w = \|g-f_c\|^2_w + \|f_c\|^2_w \geq
\|f_c\|^2_w\,.\] Hence, $f_c$ is the unique minimal norm element
in $\Hilb_w(\G)$.
\end{proof}
\begin{corollary}
The minimal norm interpolation of the sequence $(f(\lam))$ for
some $f\in\Hilb_w(\G)$ sampled on a lattice $\Lam$ in $\G$
coincides with the orthogonal projection of $f$ onto
$V_\Lam(\phi)$, say $Pf$.
\end{corollary}
\begin{proof}
This follows immediately from
\[ Pf(\lam) = \inner{T_\lam\phi}{Pf} = \inner{PT_\lam\phi}{f} =
\inner{T_\lam\phi}{f} = f(\lam)\,.  \]
\end{proof}


\section{Orthogonal Projection}

In the final section we show that the minimal norm interpolation
for a finite number of lattice elements converges to the minimal
norm interpolation for $\Lam$ when increasing the number of
lattice elements.

Let $\Lam$ be a lattice of $\G$. We denote by $\{\Lam_F\}$ a
nested sequence of finite subsets of $\Lam$ with $\bigcup_F\Lam_F
= \Lam$. We define \[ V_F(\phi) = \mbox{\rm
span}\set{T_\lam\phi}{\lam\in\Lam_F}\,.\] We recall that
$\set{T_\lam\phi}{\lam\in\Lam_F}$ is obviously a Riesz basis for
$V_F(\phi)$ whose dual basis is given by $\set{P_FT_\lam\psi}
{\lam\in\Lam_F}$ where $P_F$ denotes the orthogonal projection
onto $V_F(\phi)$. It is obvious that $P_FT_\lam\psi = 0$  for all
$\lam\not\in F$.

For a fixed element $c=(c_\lam)\in\ell^2(\Lam)$ we have seen that
\[ g = \sum_{\lam\in\Lam} c_\lam T_\lam \psi \] is the minimal
norm interpolation on $\Lam$. Set $c_F=(c_\lam)_{\lam\in\Lam_F}$.
By the same arguments used in Theorem \ref{th:minp}, we easily
deduce that the minimal norm element for $c_F$ in $\Hilb_w(\G)$ is
given by
 \[
 g_F = \sum_{\lam\in\Lam_F}c_\lam P_F T_\lam \psi = P_F
 \Big(\sum_{\lam\in\Lam_F}c_\lam T_\lam \psi\Big) =  P_F
 \Big(\sum_{\lam\in\Lam}c_\lam T_\lam \psi\Big) = P_F g\,.
 \]
In other words, the minimal norm interpolation for $c_F$ is just
the orthogonal projection of the minimal norm interpolation $g$
for $c$.

Exploiting the nested structure of the subspaces, a standard
argument in wavelet theory, e.g., \cite{Wal02}, immediately
implies
\begin{equation}\label{eq:conv}
\|g - g_F\|_w \rightarrow 0 \quad \mbox{for $F$ increasing} \,.
\end{equation}
Since $\Hilb_w(\G)$ is continuously embedded in $C(\G)$, the
convergence in (\ref{eq:conv}) holds true also for the sup-norm.

\newpage

\flushleft{

}

\end{document}